\documentstyle[amssymb,amsfonts]{amsart}

\def\be#1{\begin{equation} \label{#1}}
\def\bi{\begin{itemize}}
\def\bs{\begin{split}}
\def\es{\end{split}}
\def\ba{\begin{align}}
\def\bas{\begin{align*}}
\def\ea{\end{align}}
\def\eas{\end{align*}}
\def\R{{\hbox{\bf R}}}

\def\R{{\hbox{\bf R}}}

\def\Z{{\hbox{\bf Z}}}

\def\eps{\varepsilon}

\newenvironment{proof}{\noindent {\bf Proof} }{\endprf\par}
\def \endprf{\hfill  {\vrule height6pt width6pt depth0pt}\medskip}
\def\emph#1{{\it #1}}
\def\textbf#1{{\bf #1}}

\theoremstyle{plain}
  \newtheorem{theorem}[subsection]{Theorem}
  
  \newtheorem{proposition}[subsection]{Proposition}
  \newtheorem{lemma}[subsection]{Lemma}

\theoremstyle{remark}

\theoremstyle{definition}
  \newtheorem{definition}[subsection]{Definition}

\include{psfig}

\begin{document}

\title[Global regularity of wave maps I]{Global regularity of wave maps I.  Small critical Sobolev norm in high dimension}
\author{Terence Tao}
\address{Department of Mathematics, UCLA, Los Angeles CA 90095-1555}
\email{ tao@@math.ucla.edu}
\subjclass{35J10}

\vspace{-0.3in}
\begin{abstract}
We show that wave maps from Minkowski space $\R^{1+n}$ to a sphere $S^{m-1}$ are globally smooth if the initial data is smooth and has small norm in the critical Sobolev space $\dot H^{n/2}$, in the high-dimensional case $n \geq 5$.  A major difficulty, not present in earlier results in this area, is that the $\dot H^{n/2}$ norm barely fails to control $L^\infty$, potentially causing a logarithmic divergence in the non-linearity; however this can be overcome by using co-ordinate frames adapted to the wave map by approximate parallel transport.  In the sequel \cite{tao:wavemap2} of this paper we address the more interesting two-dimensional case $n=2$, which is energy-critical.
\end{abstract}

\maketitle

\section{Introduction}

Throughout this paper $m \geq 2$, $n \geq 1$ will be fixed integers, and all constants may depend on $m$ and $n$.

Let $\R^{1+n}$ be $n+1$ dimensional  Minkowski space  
with flat metric $\eta := \text{diag}(-1, 1, \ldots, 1)$, and let $S^{m-1} \subset \R^m$ denote the unit sphere in the Euclidean space $\R^m$.  Elements $\phi$ of $\R^m$ will be viewed as column vectors, while their adjoints $\phi^\dagger$ are row vectors.  We let $\partial_\alpha$ and $\partial^\alpha$ for $\alpha = 0,\ldots,n$ be the usual derivatives with respect to the Minkowski metric $\eta$.  We let $\Box := \partial_\alpha \partial^\alpha = \Delta - \partial_t^2$ denote the D'Lambertian.  We shall also use $\dot \phi$ for $\partial_t \phi$.

Define a \emph{wave map} to be any function $\phi$ defined on an open set in $\R^{1+n}$ taking values on the sphere $S^{m-1}$ which obeys the equation
\be{wavemap-eq}
\partial_\alpha \partial^\alpha \phi = - \phi \partial_\alpha \phi^\dagger \partial^\alpha \phi
\end{equation}
in the sense of distributions.  In order to make sense of \eqref{wavemap-eq} we shall require $\phi$ to be in $C^1_t L^2_x \cap C^0_t H^1_x$; in our high-dimensional context this regularity shall be easy to obtain.

We shall define a \emph{classical wave map} to be a wave map which is smooth and equal to a constant outside of a finite union of light cones.

For any time $t$, we use $\phi[t] := (\phi(t), \dot \phi(t))$ to denote the position and velocity of $\phi$ at time $t$.  We refer to $\phi[0]$ as the \emph{initial data} of $\phi$.  We shall always assume that the initial data $\phi[0] = (\phi(0), \dot \phi(0))$ satisfies the consistency conditions
\be{consist}
\phi^\dagger(0) \phi(0) = 1; \quad \phi^\dagger(0) \dot \phi(0) = 0
\end{equation}
(i.e. $\phi[0]$ lies on the sphere).  It is easy to show (e.g. by Gronwall's inequality) that this consistency condition is maintained in time, for smooth solutions at least.

Let $H^s := (1 + \sqrt{-\Delta})^{-s} L^2(\R^n)$ denote the usual\footnote{Strictly speaking, one cannot use $H^s$ spaces for functions on the sphere, since they cannot globally be in $L^2$.  To get around this we shall abuse notation and allow constant functions to lie in $H^s$ with zero norm whenever the context is for functions on the sphere.  Thus when we say that $\phi(t)$ is in $H^s$, we really mean that $\phi(t) - c$ is in $H^s$ for some constant $c$.}  $L^2$ Sobolev spaces.  Since the equation \eqref{wavemap-eq} is invariant under the scaling $\phi \mapsto \phi_\lambda$ defined by
$$ \phi_\lambda(t,x) := \phi(t/\lambda, x/\lambda)$$
we see that the critical regularity is $s=n/2$.

The Cauchy problem for wave maps has been extensively studied (see references); we refer the interested reader to the surveys in \cite{kman.barrett}, \cite{kman.selberg:survey}, \cite{shatah-struwe}, \cite{struwe.barrett}.   For sub-critical regularities $s > n/2$ it is known (see \cite{kl-mac:null3}, \cite{kman.selberg}, \cite{keel:wavemap} for the $n \geq 4$, $n=2,3$, and $n=1$ cases respectively) that the Cauchy problem for \eqref{wavemap-eq} is locally well-posed for initial data $\phi[0]$ in $H^s \times H^{s-1}$, and the solution can be continued (without losing regularity) as long as the $H^s$ norm remains bounded.   The critical result however is more subtle.  Well-posedness and regularity was demonstrated in the critical Besov space $\dot B^{n/2}_1$ in \cite{tataru:wave1} in the high-dimensional case $n \geq 4$ and in \cite{tataru:wave2} for $n=2,3$; in the one-dimensional case $n=1$, a logarithmic cascade from high frequencies to low frequencies causes ill-posedness in the critical Besov and Sobolev spaces \cite{tao:ill}, although global smooth solutions can still be constructed thanks to the sub-critical nature of the energy norm (see e.g. \cite{shatah.zurich}).  As is to be expected at the critical regularity, these results give a global well-posedness (and regularity and scattering) when the norm of the initial data is small.

The question still remains as to whether the wave map equation \eqref{wavemap-eq} is well-posed in the critical Sobolev space $\dot H^{n/2} := \sqrt{-\Delta}^{-n/2} L^2$ in two and higher dimensions, with global well-posedness and regularity expected for small data.  This question is especially interesting in the two dimensional case, since the critical Sobolev space is then the energy norm $H^1$, and one also expects to exploit conservation of energy (and some sort of energy non-concentration result) to obtain global well-posedness and regularity for data with large energy.  (In dimensions three and higher one does not have large data global well-posedness for the sphere, even for smooth symmetric data; see \cite{shatah.shadi.blow}). However the Sobolev space $\dot H^{n/2}$ is significantly less tractable than its Besov counterpart $\dot B^{n/2}_1$; for instance, $\dot H^{n/2}$ norm fails to control the $L^\infty$ norm due to a logarithmic pile-up of frequencies.  This logarithmic divergence is responsible for failure of any strengthened version of well-posedness (uniform, Lipschitz, or analytic) for wave maps at this regularity, as well as ill-posedness in very similar equations, and is a serious obstacle to any iteration-based argument.  See \cite{kman.selberg}, \cite{keel:wavemap}, \cite{nak:wavemap} for further discussion.

Our main result is the following.

\begin{theorem}\label{main2}
Let $n \geq 5$ and $s > n/2$, and suppose that the initial data $\phi[0]$ is in $H^s \times H^{s-1}$ and has sufficiently small $\dot H^{n/2} \times \dot H^{n/2-1}$ norm.  Then the solution to the Cauchy problem for \eqref{wavemap-eq} with initial data $\phi[0]$ can be continued in $H^s \times H^{s-1}$ globally in time.  In particular, smooth solutions stay smooth when the initial data has small $\dot H^{n/2} \times \dot H^{n/2-1}$ norm.  Furthermore, if $|s - n/2| < 1/2$, we have the global bounds\footnote{We of course adopt the convention that $A \lesssim B$ denotes the inequality $A \leq CB$ for some constant $C$ depending only on $n$, $m$.}
\be{bound}
\| \phi[t] \|_{L^\infty_t (\dot H^s_x \times \dot H^{s-1}_x)} \lesssim \| \phi[0] \|_{\dot H^s_x \times \dot H^{s-1}_x}.
\end{equation}
\end{theorem}

Our argument also shows that $\phi$ obeys the expected range of Strichartz estimates globally in spacetime, although we will not write down a precise statement here.

Our arguments are heavily based on the geometric structure of the equation \eqref{wavemap-eq}; in particular, they do not directly apply to the associated difference equation.  As a consequence we have not been able to obtain a well-posedness\footnote{We should remark at this point that strong versions of well-posedness, such as uniform, Lipschitz, or analytic well-posedness, are known to fail at the critical Sobolev regularity (see \cite{keel:wavemap}, \cite{nak:wavemap}).  Presumably one would have to renormalize the difference equation in a manner strongly dependent on the initial data.}  result at the critical regularity $\dot H^{n/2} \times \dot H^{n/2-1}$, even for small data. The argument also does not directly yield a scattering result, although this obstruction seems to be less serious.  We will not pursue these matters.

Our main tools are Littlewood-Paley decomposition and Strichartz estimates (as in \cite{tataru:wave1}), combined with some geometric identities (such as $\phi^\dagger \partial_\alpha \phi = 0$) and the use of a good co-ordinate frame constructed by approximate parallel transport; this renormalization is crucial in order to remove the possibility of logarithmic divergence.  The high-dimensional case $n \geq 5$ is significantly easier than the low-dimensional cases, because of the strong decay of the wave equation ($t^{-2}$ or better) as well as the rarity of parallel interactions (one expects interactions of angle $O(\theta)$ to only occur $O(\theta^{n-1}) = O(\theta^4)$ of the time).  Because of these advantages, we shall not need to develop bilinear estimates to take advantage of the null form structure of \eqref{wavemap-eq} (cf. \cite{tat:5+1}), although we shall heavily exploit the geometric structure of this equation.  Instead of bilinear estimates, we shall rely primarily on Strichartz estimates, and in particular on the $L^2_t L^\infty_x$-type and $L^2_t L^4_x$-type estimates which are not available in low dimensions.  With such strong estimates it shall be easy (after applying the renormalization) to obtain $L^1_t L^2_x$-type estimates on the non-linearity, so that one can then close the argument by energy estimates.  (Without the normalization, Strichartz estimates only work when one has half a derivative more than critical).

In the $n=4$ case one loses the $L^2_t L^4_x$ estimate, however one can introduce $X^{s,b}$ spaces to obtain $L^2_t L^2_x$ estimates on $\Box \phi$.  This, together with the identity
$$ 2\partial_\alpha \phi \partial^\alpha \psi = \Box(\phi \psi) - \phi \Box(\psi) - \Box(\phi) \psi$$
to exploit the null structure, should be able to cover the $n=4$ case in analogy with the arguments in \cite{tataru:wave1}, however we have elected not to do this to keep the argument as simple as possible.  

In the $n=2,3$ cases one also loses the $L^2_t L^\infty_x$ estimate, which seems to defeat any attempt to prove these results purely by standard Strichartz estimates (although the $n=3$ case is probably salvageable for radial data, thanks to the endpoint Strichartz estimate holding in that context).  Fortunately, the low dimensional case has been effectively handled (in the Besov space case) by the more sophisticated arguments of \cite{tataru:wave2}, using the additional ingredients of angular frequency decomposition and $X^{s,b}$-type spaces (to more effectively exploit the null form structure in \eqref{wavemap-eq}) and null frames (to recover estimates of $L^2_t L^\infty_x$ type).  In the sequel \cite{tao:wavemap2} to this paper we shall adapt the arguments in \cite{tataru:wave2} to cover the $n=2$ small energy case, as well as the remaining cases $n=3,4$.

The main novel ingredient in our approach is the use of adapted co-ordinate frames constructed by approximate parallel transport along (Littlewood-Paley regularizations of) $\phi$.  The construction presented here is heavily dependent of the geometry of the sphere, although this should in principle extend to other compact manifolds by using the machinery of Helein in his work \cite{helein} on harmonic maps\footnote{Since the preparation of this manuscript, the author has learnt (Klainerman, personal communication) that the arguments here have been successfully extended to arbitrary Lie groups by Klainerman and Rodnianski.}. 

Without the use of these frames, the usual iteration approach for \eqref{wavemap-eq} fails at the critical regularity because of a logarithmic pile-up of high-low frequency interactions.  The effect of the adapted co-ordinate frame is to transform the high-low frequency interaction into other terms which are more tractable, such as high-high frequency interactions, or high-low interactions in which a derivative has been moved from a high-frequency term to a low-frequency one.

In the remainder of this section we shall informally motivate the key ideas in the argument.  In doing so we shall make frequent use of the following heuristic: if $\phi$, $\psi$ are two functions, and $\psi$ is much rougher (i.e. higher frequency) than $\phi$, then $(\nabla \phi) \psi$ is very small compared to $\phi \nabla \psi$.  In other words, we should be able to neglect terms in which derivatives fail to fall on rough functions, and land instead on smooth ones.  (Indeed, these terms can usually be treated just by Strichartz estimates).  In particular, we expect to have $\nabla(\phi \psi) \approx \phi \nabla \psi$ (which can be viewed as a statement that $\phi$ is approximately constant when compared against $\psi$).

Let us suppose that our wave map $\phi$ has the form $\phi = \tilde \phi + \eps \psi$, where $\tilde \phi$ is a smooth wave map, $0 < \eps \ll 1$ and $\psi$ is a $H^{n/2}$ function which is much rougher than $\tilde \phi$.  (In other words, $\phi$ is a small rough perturbation of a smooth wave map).  If we ignore terms which are quadratic or better in $\eps$, or which fail to differentiate the rough function $\psi$, we obtain the linearized equation
\be{linearized}
\partial_\alpha \partial^\alpha \psi = - 2\tilde \phi \partial_\alpha \tilde \phi^\dagger \partial^\alpha \psi
\end{equation} 
for $\psi$.  Also, since $\tilde \phi$ and $\tilde \phi + \eps \psi$ both take values on the sphere we see that 
\be{nabla}
\tilde \phi^\dagger \psi = 0; \quad \tilde \phi^\dagger \partial_\alpha \psi = 0
\end{equation}
(again ignoring terms quadratic in $\eps$, and terms where the derivative fails to land on $\psi$).

In order to keep the $H^s$ norm of $\tilde \phi + \eps \psi$ from blowing up, we need to prevent the $\dot H^{n/2}$ norm from being transferred from $\tilde \phi$ to $\eps \psi$.  In particular, we need $L^\infty_t \dot H^{n/2}_x$ bounds on $\psi$ which are independent of $\eps$.  We would also like the corresponding Strichartz estimates for $\psi$, in order to control the error terms that we have been ignoring.  (This scheme is not restricted to rough perturbations of smooth wave maps, and will be adapted to general wave maps by use of Littlewood-Paley projections).

Despite being linear, the equation \eqref{linearized} is not very well-behaved, having no obvious cancellation structure (beyond the null form, which is not particularly useful in the high-dimensional setting).  In order to iterate away the first-order terms on the right-hand side of \eqref{linearized} we would like $\tilde \phi \partial_\alpha \tilde \phi^\dagger$ to be in $L^1_t L^\infty_x$.  In principle this might be feasible if we had the Strichartz estimate $\nabla^{1/2} \tilde \phi \in L^2_t L^\infty_x$, but this estimate just barely fails to hold because of a logarithmic divergence in the frequencies.  However, if we could somehow ensure that the derivative in $\tilde \phi \partial_\alpha \tilde \phi^\dagger$ always fell on a low-frequency component of $\tilde \phi$ and not on a high-frequency component then one would have a chance of iterating away the non-linearity\footnote{The author thanks Chris Sogge for this observation.}.  This will be accomplished by a renormalization using a co-ordinate frame adapted to $\phi$.  

We begin by taking advantage of \eqref{nabla} to rewrite \eqref{linearized} in a form reminiscent of parallel transport:
\be{cancel}
\partial_\alpha \partial^\alpha \psi = 2 A_\alpha \partial^\alpha \psi
\end{equation}
where $A_\alpha$ is the matrix
$$ A_\alpha := \partial_\alpha \tilde \phi \tilde \phi^\dagger - \tilde \phi \partial_\alpha \tilde \phi^\dagger.$$
Note that \eqref{cancel} exhibits more cancellation than \eqref{linearized}, as $A_\alpha$ is now anti-symmetric.  This type of trick is standard in the study of wave and harmonic maps, see e.g. \cite{keel:wavemap}, \cite{helein}, \cite{changwangyang}, etc.

To solve \eqref{cancel}, let us first consider the ODE analogue
\be{cancel-ode}
\ddot \psi = 2 A_0 \dot \psi.
\end{equation}
The matrix $A_0$ is anti-symmetric.  Thus if we let $U(t)$ be the matrix-valued function solving the ODE
$$ \dot U(t) = A_0 U(t)$$
with $U(0)$ initialized to the identity matrix (say), then we see that $\frac{d}{dt} (U U^\dagger) = 0$ and thus that $U$ remains orthogonal for all time.  Indeed, one can view $U$ as the parallel transport of the identity matrix along the trajectory of $\tilde \phi$.  Furthermore, since $\tilde \phi$ is smooth, we see that $U$ is also smooth, and in particular is much smoother than $\psi$.  One can then use the linear change of variables $\psi = Uw$, and ignore terms which fail to differentiate the rough function $w$, to rewrite \eqref{cancel-ode} as the trivial equation $\ddot w = 0$.  

The ODE example of \eqref{cancel-ode} suggests that \eqref{cancel} might be simplified by applying some orthogonal matrix $U$ to the wave $\psi$, or in other words by viewing $\psi$ in a carefully chosen co-ordinate frame.  (This fits well with the corresponding experience of harmonic maps in \cite{helein}).  Ideally, we would like $U$ to be carried by parallel transport by $\tilde \phi$ in all directions.  More precisely, we would like $U$ to solve the PDE
\be{ook}
\partial_\alpha U = A_\alpha U
\end{equation}
for each $\alpha$. If we make the improbable assumption that $U$ obeyed \eqref{ook} exactly for all $\alpha$, we can then substitute $\psi = Uw$ as before and ignore all terms which fail to differentiate the rough function $w$ to transform \eqref{cancel} to the free wave equation
$$
\partial_\alpha \partial^\alpha w = 0$$
which we of course know how to solve.  

Unfortunately, the system \eqref{ook} of PDE is overdetermined, and in general has no solution (since the parallel transport connection induced by $\tilde \phi$ will have a small\footnote{More precisely, the curvature only contains terms which are quadratic in the first derivatives of $\tilde \phi$, as opposed to being linear in the second derivatives of $\tilde \phi$.  This phenomenon seems specific to the wave maps equation; if one tries to apply the techniques here to (for instance) the Maxwell-Klein-Gordon or Yang-Mills equations at the critical Sobolev regularity, an obstruction arises because the connection $A$ has no reason to have a good curvature, regardless of the choice of gauge.  At best one can place these equations in the Coulomb gauge, which was already known to be the most useful gauge to study these equations.} but non-zero curvature).  Nevertheless, it is possible to use Littlewood-Paley  theory to construct a satisfactory \emph{approximate} solution $U$ to \eqref{ook}.  Specifically, we perform the Littlewood-Paley  decomposition $\tilde \phi = \phi_{-M} + \sum_{-M < k} \phi_k$, where $M$ is a large number, $\phi_{-M}$ is the portion of $\phi$ on frequencies $|\xi| \lesssim 2^{-M}$, and $\phi_k$ is the portion on frequencies $|\xi| \sim 2^k$.  We then define $U = U_{-M} + \sum_{-M < k} U_k$, where $U_{-M}$ is the identity matrix, and the $U_k$ are defined recursively by the formula
\be{ack}
U_k := (\phi_k \phi_{<k}^\dagger - \phi_{<k} \phi_k^\dagger) U_{<k}
\end{equation}
where $\phi_{<k}$, $U_{<k}$ are the functions
$$ \phi_{<k} := \phi_{-M} + \sum_{-M < k' < k} \phi_k, \quad U_{<k} := U_{-M} + \sum_{-M < k' < k} U_k.$$
It then transpires that the matrix $U$ is approximately orthogonal and approximately satisfies \eqref{ook}, provided that the $\dot H^{n/2}$ norm of $U$ is sufficiently small and $M$ is sufficiently large.  The point is that $\phi_k$ is a rougher function than $\phi_{<k}$, and so one can (heuristically) neglect terms where the derivative falls on $\phi_{<k}$ instead of $\phi_k$.  Similarly for $U_k$ and $U_{<k}$.  Thus we can morally differentiate \eqref{ack} to obtain
\be{ack-diff}
\partial_\alpha U_k \approx (\partial_\alpha \phi_k \phi_{<k}^\dagger - \phi_{<k} \partial_\alpha \phi_k^\dagger) U_{<k}
\end{equation}
and \eqref{ook} follows by summing the telescoping series (and continuing to neglect the same type of terms as before).  The approximate orthogonality of $U$ is based on the observation (from \eqref{ack}) that $U_k^\dagger U_{<k}+  U_{<k}^\dagger U_k = 0$.  Summing this in $k$ and telescoping, we obtain
$$ U^\dagger U = I + \sum_{k > -M} U^\dagger_k U_k.$$
The summation on the right-hand side then turns out to be negligible if we assume $\tilde \phi$ is small in $\dot H^{n/2}$, since this implies from Sobolev embedding that the $L^\infty$ norms of the $\phi_k$ (and hence the $U_k$) are small in $l^2$.  (A similar argument can be used to dispose of the error terms which were neglected in \eqref{ack-diff}).
If one then transforms \eqref{cancel} using $\psi = Uw$ as before, we obtain a non-linear wave equation for $w$, but all the terms in the non-linearity either contain expressions such as $\sum_k U_k U_k^\dagger$ which are quadratic\footnote{Basically, such quadratic expressions effectively improve the Sobolev space $\dot H^{n/2}$ to the Besov space $\dot B^{n/2}_1$, which in principle can be treated by the arguments in \cite{tataru:wave1}.} in the frequency parameter $k$, or have all derivatives falling on smooth functions rather than rough ones.  Both types of terms turn out to be easily controlled by Strichartz estimates.

This work was conducted at UCLA, Tohoku University, UNSW, and the French Alps.  The author thanks Daniel Tataru, Mark Keel, and Sergiu Klainerman for very helpful discussions, insights, and encouragement, and to Sergiu Klainerman, Kenji Nakanishi, Igor Rodnianski, and the referee for pointing out errors in an early preprint.  The author is a Clay Prize Fellow and is supported by grants from the Sloan and Packard foundations.
 
\section{Littlewood-Paley projections and Strichartz estimates}\label{strichartz-sec}

In this section we set out notation for two basic tools in this argument.

We begin with Littlewood-Paley operators.  If $\phi(t,x)$ is a function in spacetime, we define the spatial Fourier transform $\hat \phi(t,\xi)$ by
$$ \hat \phi(t,\xi) := \int_{\R^n} e^{-2\pi i x \cdot \xi} \phi(t,x)\ dx.$$
Fix $m(\xi)$ to be a non-negative radial bump function supported on $|\xi| \leq 2$ which equals 1 on the ball $|\xi| \leq 1$.  For each integer $k$, we define the Littlewood-Paley projection operators $P_{\leq k} = P_{<k+1}$ to the frequency ball $|\xi| \lesssim 2^k$ by the formula
$$ \widehat {P_{\leq k} \phi}(t,\xi) := m(2^{-k}\xi) \hat \phi(t,\xi),$$
and the projection operators $P_k$ to the frequency annulus $|\xi| \sim 2^k$ by the formula
$$ P_k := P_{\leq k} - P_{<k}.$$

We also define more general projections $P_{k_1 \leq \cdot \leq k_2}$ by
$$ P_{k_1 \leq \cdot \leq k_2} := P_{\leq k_2} - P_{< k_1}.$$
Similarly define $P_{k_1 < \cdot \leq k_2}$, etc.

Note that if $\phi$ is a smooth function which is equal to a constant $e$ outside of a compact set, then we have the Littlewood-Paley decomposition
\be{lp}
\phi = e + \sum_k P_k \phi.
\end{equation}
Also, we remark that the Littlewood-Paley  projections defined above commute with all constant-coefficient differential operators and are bounded on every Lebesgue space (including mixed-norm spacetime Lebesgue spaces).

Because we are in the high-dimensional case $n \geq 5$, we will not need $X^{s,b}$-type spaces.  Indeed, these spaces do not quite seem to be the right tool for dealing with critical Sobolev regularity problems (despite being very powerful for subcritical problems).  Because we are avoiding these spaces, we may localize in time freely without encountering distracting technicalities involving the temporal Fourier transform.

We now describe the Strichartz estimates that we need.
Let us call a pair $(q,r)$ of exponents \emph{admissible} if $2 \leq q, r \leq \infty$ and 
$$ \frac{1}{q} + \frac{(n-1)/2}{r} \leq \frac{(n-1)/2}{2}.$$
For any integer $k$, we define the \emph{($\dot H^{n/2}$-normalized) Strichartz space at frequency $2^k$}, $S_k(\R^{1+n})$, to be the space of functions on spacetime whose norm is given by\footnote{The powers of $2^k$ which will appear in the sequel are not mysterious, and can be explained by scaling. One should think of $2^k$ has having the units of frequency (i.e. inverse length), so that $S_k$ has the scaling of $\dot H^{n/2}_x$ or $L^\infty_t L^\infty_x$.} 
$$ \| \phi \|_{S_k} := \sup_{q,r} 2^{\frac{k}{q} + \frac{kn}{r}} (\| \phi \|_{L^q_t L^r_x} + 2^{-k} \| \partial_t \phi \|_{L^q_t L^r_x}),$$
where the supremum ranges over all admissible exponents $(q,r)$.  Similarly define $S_k(I \times \R^n)$ for time intervals $I$.  Generally speaking, the large values of $r$ are good for low-frequency terms, and conversely for high-frequency terms.  In our high dimensional setting $n \geq 5$ we have a very large set of Strichartz estimates which will be more than adequate for our purposes.

We shall only use specific values of $q$ and $r$ in our argument.  More precisely, we observe that control of the $S_k$ norm gives the estimates
\begin{align}
\| \phi \|_{L^2_t L^{2(n-1)/(n-3)}_x} &\leq 2^{-\frac{k}{2} - \frac{nk}{2} + 
\frac{nk}{n-1}} \| \phi \|_{S_k} \label{endpoint}\\
\| \phi \|_{L^2_t L^4_x} &\leq 2^{-\frac{k}{2} - \frac{nk}{4}} \| \phi \|_{S_k} \label{four}\\
\| \phi \|_{L^2_t L^{n-1}_x} &\leq 2^{-\frac{k}{2} - \frac{nk}{n-1}} \| \phi \|_{S_k} \label{dual}\\
\| \phi \|_{L^2_t L^\infty_x} &\leq 2^{-\frac{k}{2}} \| \phi \|_{S_k}; \label{2infty}\\
\| \phi \|_{L^4_t L^{2(n-1)}_x} &\leq 2^{-\frac{k}{4} - \frac{nk}{2(n-1)}} \| \phi \|_{S_k} \label{4dual}\\
\| \phi \|_{L^\infty_t L^\infty_x} &\leq \| \phi \|_{S_k} \label{infty}\\
\| \phi \|_{L^\infty_t L^2_x} &\leq 2^{-\frac{nk}{2}} \| \phi \|_{S_k} \label{energy}
\end{align}
note that we have used the hypothesis $n \geq 5$ in order to obtain the admissibility of \eqref{four} and \eqref{dual}.  Clearly, one can also estimate the time derivative $\phi_t$ in the above norms by paying an additional power of $2^k$.

Suppose $\phi$ is a smooth function and $\psi$ is a rough function.  Our non-linearity is cubic in $\phi$ and $\psi$ with two derivatives somewhere.  To estimate a term such as $\phi \nabla \phi \nabla \psi$ in $L^1_t L^2_x$, we shall usually estimate $\nabla \psi$ using \eqref{energy} and $\phi$, $\nabla \phi$ using \eqref{2infty}; this turns out to work as long as the $\phi$ term has equal or higher frequency to $\nabla \phi$.  The estimate \eqref{four} is useful for obtaining $L^1_t L^2_x$ control on terms such as $\nabla \psi \nabla \psi \phi$ which are quadratic in the high frequencies, and the pair \eqref{dual}, \eqref{endpoint} are useful for controlling terms such as $\phi \nabla^2 \phi \psi$.  Finally, the triplet \eqref{4dual}, \eqref{4dual}, \eqref{endpoint} can handle the term $\nabla \phi \nabla \phi \psi$.  There is a certain amount of flexibility in our choice of exponents, especially in large dimensions; in particular, the endpoint \eqref{endpoint} can be avoided when $n \geq 6$.

We have the Strichartz estimates (see e.g. \cite{tao:keel} and the references therein):

\begin{theorem}[Strichartz estimates]\label{strich} Let $k$ be an integer.  For any function $\phi$ on $\R^{1+n}$ with Fourier support on $|\xi| \sim 2^k$, we have
$$ \| \phi \|_{S_k} \lesssim \| \phi[0] \|_{\dot H^{n/2}_x \times \dot H^{n/2-1}_x} + 2^{-k+\frac{nk}{2}} \| \Box \phi \|_{L^1_t L^2_x}.$$
Similarly if $\R^{1+n}$ is replaced by $I \times \R^n$ for any interval $I$ containing the origin. 
\end{theorem}

One could place $\Box \phi$ in other spaces than $L^1_t L^2_x$ (indeed, one could use the dual of $S_k$, appropriately normalized) but we shall not need to do so here.

\section{The main proposition}\label{main-sec}

In this section we state the main proposition which will be used to prove Theorem \ref{main2}.  Roughly speaking, the proposition asserts that if the $\dot H^{n/2}$ norm of a (classical) wave map is initially small, then it stays small for all time, and its frequency profile does not change substantially.  From this proposition and the existing well-posedness theory it shall be an easy matter to obtain Theorem \ref{main2} by general arguments.

Throughout the paper we fix $\sigma$ to be a constant depending only on $n$ such that $0 < \sigma < 1/2$ (e.g. $\sigma := 1/4$ will do).  We also fix $0 < \eps \ll 1$ to be a small constant depending only on $n$, $m$, $\sigma$ ($\eps := 2^{-100nm}$ will suffice).

\begin{definition}\label{envelope-def}  A \emph{frequency envelope} is a sequence $c = \{c_k\}_{k \in \Z}$ of positive reals such that we have the $l^2$ bound
\be{l2-ass}
\|c\|_{l^2} \lesssim \eps
\end{equation}
and the local constancy condition
\be{local}
2^{-\sigma |k-k'|} c_{k'} \lesssim c_k \lesssim 2^{\sigma |k-k'|} c_{k'}
\end{equation}
for all $k, k' \in \Z$.  In particular we have $c_k \sim c_k'$ whenever $k = k' + O(1)$.  If $c$ is a frequency envelope and $(f,g)$ is a pair of functions on $\R^n$, we say that $(f,g)$ \emph{lies underneath the envelope} $c$ if one has
$$ \| P_k f \|_{\dot H^{n/2}} + \| P_k g \|_{\dot H^{n/2-1}} \leq c_k$$
for all $k \in \Z$.
\end{definition}

Note that if $(f,g)$ lies underneath an envelope $c$, then 
\be{init}
\| (f,g) \|_{\dot H^{n/2} \times \dot H^{n/2-1}} \lesssim \eps.
\end{equation}
Conversely, if \eqref{init} holds (with slightly better implicit constants) then there exists an envelope $c$ which $f$ lies underneath, for instance one can take
\be{formula}
c_k := \sum_{k' \in \Z} 2^{-\sigma|k-k'|} (\| P_{k'} f \|_{\dot H^{n/2}} + \| P_{k'} g \|_{\dot H^{n/2-1}}).
\end{equation}

\begin{proposition}[Main Proposition]\label{main-prop}  Let $0 < T < \infty$, $c$ be a frequency envelope, and $\phi$ be a classical wave map on $[0,T] \times \R^n$ such that $\phi[0]$ lies underneath the envelope $c$.  Then, if $\eps$ is sufficiently small, we have the bounds
\be{control}
\| P_k \phi \|_{S_k([0,T] \times \R^n)} \leq C_0 c_k
\end{equation}
for all $k \in \Z$, where $C_0 \gg 1$ is an absolute constant depending only on $n$, $m$ (i.e. independent of $T$, $\phi$, $c$, $\eps$).  From \eqref{control}, \eqref{energy}, and Definition \ref{envelope-def} we have in particular that $\phi[t]$ lies underneath the envelope $C_0 c$ for all $t \in [0,T]$.
\end{proposition}

We remark that the finiteness of $T$ is needed for some continuity and limiting arguments to work, but otherwise $T$ plays no role in the estimates.  At first glance this Proposition seems to be merely a variant of a conservation law for the $\dot H^{n/2}$ norm, but the fact that the envelope $c$ is arbitrary makes this result far more powerful; indeed, it effectively allows one to deduce the $\dot H^{n/2}$ regularity theory from the $H^{n/2+}$ theory.  Also, we remark that if we strengthen the $l^2$ control in \eqref{l2-ass} to $l^1$ control, then Proposition \ref{main-prop} was essentially proven in \cite{tataru:wave1}.

Assuming this Proposition, the proof of Theorem \ref{main2} is now easy.  In light of the existing regularity and well-posedness results\footnote{Strictly speaking, this regularity result was stated for data with small $H^s \times H^{s-1}$ norm, but the result extends as well to large data; see \cite{kman.selberg:survey}, \cite{selberg:thesis}.  Alternatively, one can exploit finite speed of propagation and the fact (which is a variant of Poincare's inequality) that if $\phi$ has large $H^s$ norm, then the $H^{s'}$ norm of $\phi$ on a small ball $B(x,r)$ decays like a power of $r$ for $n/2 < s' < s$, recalling that we allow a constant to be subtracted from $\phi$ when defining the $H^{s'}$ norm.} in \cite{kman.selberg} (see also \cite{tataru:wave1}), as well as finite speed of propagation, it suffices to prove that the theorem when $s$ is close to $n/2$, or specifically when $n/2 < s < n/2+\sigma$, and we may also assume by standard limiting arguments that $\phi$ is a classical wave map.  But if the initial data $\phi[0]$ is in $H^s$, then there exists an envelope $c$ that $\phi[0]$ lies under which decays like $2^{-(s-n/2)k}$ as $k \to \infty$ (just by applying \eqref{formula}).  By Proposition \ref{main-prop}, $\phi[t]$ lies under the envelope $Cc_0$ for as long as the solution stays smooth, which in turn implies that the $H^s$ norm of $\phi[t]$ stays uniformly bounded for as long as the solution stays smooth.  The claim then follows from the regularity results in \cite{kman.selberg}, \cite{tataru:wave1}.  The bound \eqref{bound} follows easily as a by-product of the above argument (note that this bound also applies for $s$ slightly below $n/2$).

We close this section by performing some preliminary reductions for Proposition \ref{main-prop}.  Fix $\phi$; we may assume that $\phi(0) = e_1$ and $\partial_t \phi(0) = 0$ outside of a compact set.  Since $\phi$ is a classical wave map, the $\dot H^{n/2} \times \dot H^{n/2-1}$ norm of $P_k \phi$ decays faster than $2^{-|k|}$ as $k \to \pm \infty$.  We may thus assume that $c$ decays like $2^{-\sigma |k|}$ as $k \to \infty$.

Fix $c$.  We shall apply the continuity method. Let $A \subseteq [0,T]$ be the set of all times $T' \subset [0,T]$ such that \eqref{control} holds (with $T$ replaced by $T'$).  Observe that $A$ is closed.  Let $A' \subseteq [0,1]$ be defined similarly but with \eqref{control} replaced by the variant
\be{control-2}
\| P_k \phi \|_{S_k([0,T] \times \R^n)} \leq 2C_0 c_k
\end{equation}
and with $T$ replaced again by $T'$.  Clearly $A \subset A'$.  In fact, since $\phi$ is a classical wave map we see that $\| P_k \phi \|_{S_k}$ decays like $O(2^{-|k|})$ as $k \to \pm \infty$, and we can therefore make the stronger statement that $A$ is contained in the interior of $A'$.  Also, $A$ trivially contains 0.  Thus if we can prove that $A' \subseteq A$, then $A$ would be open, closed, and non-empty in $[0,T]$, so that $T \in A$.  This will prove Proposition \ref{main-prop}.

From the above discussion, we have reduced Proposition \ref{main-prop} (and Theorem \ref{main2}) to

\begin{proposition}[Reduced Main Proposition]\label{reduced}  Let $c$ be a frequency envelope, $0 < T < \infty$, and let $\phi$ be a classical wave map on $[0,T] \times \R^n$ such that $\phi[0]$ lies underneath $c$, and that \eqref{control-2} holds for all $k$.  Then \eqref{control} holds for all $k$ (assuming that $C_0$ is sufficiently large and $\eps$ sufficiently small).
\end{proposition}

The proof of Proposition \ref{reduced} shall occupy the next few sections.  In Section \ref{linear-sec} we apply Littlewood-Paley projections to reduce \eqref{wavemap-eq} to the linearized variant \eqref{linearized}, modulo acceptable errors.  In Section \ref{tangent-sec} we develop rigorous versions of \eqref{nabla}, and thus convert \eqref{linearized} to the anti-symmetric form \eqref{cancel}.  In Section \ref{transport-sec} we use the approximate parallel transport scheme \eqref{ack} to construct a good co-ordinate frame $U$, and in Section \ref{conclusion-sec} we show how the transformation $\psi = Uw$ transforms \eqref{cancel} into the free wave equation, modulo acceptable errors.

\section{Linearization}\label{linear-sec}

Fix $c$, $T$, $\phi$, and suppose that the hypotheses of Proposition \ref{reduced} hold.  Henceforth all spacetime norms will be assumed to be on the slab $[0,T] \times \R^n$.

Since the differential operator $2^{-k} \nabla$ is bounded on frequencies $|\xi| \sim 2^k$, we thus have that
\be{deriv}
\| \nabla^j P_k \phi \|_{S_k} \lesssim 2^{jk} C_0 c_k
\end{equation}
for all $k \in \Z$ and all $j$, with the implicit constant depending on $j$.

We need to show \eqref{control}.  By scale-invariance (scaling $T$, $c$, and $\phi$ appropriately) it suffices to show that
\be{psi-bound}
\| \psi \|_{S_0} \leq C_0 c_0,
\end{equation}
where we define
$$ \psi := P_0 \phi.$$

By applying $P_0$ to \eqref{wavemap-eq} we obtain
\be{psi-eq}
\Box \psi = - P_0(\phi \partial_\alpha \phi^\dagger \partial^\alpha \phi).
\end{equation}
We would like to convert this non-linear equation into the linearized equation \eqref{linearized}, modulo acceptable errors.  To define ``acceptable error'' we introduce

\begin{definition}\label{error-def}  A function $F$ on $[0,T] \times \R^n$ is said to be an \emph{acceptable error} if
$$ \|F\|_{L^1_t L^2_x} \lesssim C_0^3 \eps c_0,$$
and we shall write $F = error$ to denote this.
\end{definition}

Ideally we would like the entire non-linearity in \eqref{psi-eq} to be an acceptable error, as one could then use Theorem \ref{strich} to obtain \eqref{psi-bound}.  Although we cannot quite do this directly, we can show that almost all of the non-linearity is an acceptable error, and the remaining term can be renormalized by a suitable change of co-ordinates to also be acceptable error.

We shall use $\tilde \phi := P_{\leq -10} \phi$ to denote the regularization of $\phi$.  Since $\phi$ lies on the sphere, we thus have 
\be{infty-bound}
\| \tilde \phi \|_{L^\infty_t L^\infty_x} \lesssim \| \phi \|_{L^\infty_t L^\infty_x} = 1.
\end{equation}
Also, from \eqref{2infty}, \eqref{deriv}, and the triangle inequality we observe the useful bound
\be{2infty-bound}
\| \nabla \tilde \phi \|_{L^2_t L^\infty_x} \lesssim C_0 \eps.
\end{equation}

The main result of this section is then

\begin{proposition}\label{linear-prop}  We have
\be{proj}
P_0(\phi \partial_\alpha \phi^\dagger \partial^\alpha \phi) = 2 \tilde \phi \partial_\alpha \tilde \phi^\dagger \partial^\alpha \psi + error.
\end{equation}
In particular, from \eqref{psi-eq} we have
\be{proj-effect}
\Box \psi = -2 \tilde \phi \partial_\alpha \tilde \phi^\dagger \partial^\alpha \psi + error.
\end{equation}
\end{proposition}

\begin{proof}
We apply the Littlewood-Paley decomposition to split the expression inside the projection $P_0$ in \eqref{proj} as 
\begin{align}
&\sum_{k_1, k_2, k_3: \max(k_2,k_3) > 10, |k_2 - k_3| \leq 5} (P_{k_1} \phi)
(\partial_\alpha P_{k_2} \phi^\dagger) (\partial^\alpha P_{k_3} \phi) \label{dec1}\\
&+\sum_{k_1, k_2, k_3: \max(k_2,k_3) > 10, |k_2 - k_3| > 5} (P_{k_1} \phi)
(\partial_\alpha P_{k_2} \phi^\dagger) (\partial^\alpha P_{k_3} \phi) \label{dec2}\\
&+\sum_{k_1, k_2, k_3: \max(k_2,k_3) < 10, k_1 > -10} (P_{k_1} \phi)
(\partial_\alpha P_{k_2} \phi^\dagger) (\partial^\alpha P_{k_3} \phi) \label{dec3}\\
&+ \tilde \phi \partial_\alpha \tilde \phi^\dagger \partial^\alpha \tilde \phi
\label{dec4}\\
&+ \tilde \phi (\partial_\alpha P_{-10 < \cdot < 10} \phi^\dagger) (\partial^\alpha P_{-10 < \cdot < 10} \phi). \label{dec5}\\
&+ \tilde \phi \partial_\alpha \tilde \phi^\dagger (\partial^\alpha P_{-10 < \cdot < 10} \phi) \label{dec6}\\
&+ \tilde \phi (\partial_\alpha P_{-10 < \cdot < 10} \phi^\dagger) \partial^\alpha \tilde \phi. \label{dec7}
\end{align}

As it turns out, all the terms except for \eqref{dec6}, \eqref{dec7} will be of the form $error$.

We first consider the high-frequency contributions \eqref{dec1}, \eqref{dec2}.  We can rewrite \eqref{dec1} as
$$ \phi \sum_{k_2, k_3: \max(k_2,k_3) > 10, |k_2 - k_3| \leq 5} (\partial_\alpha P_{k_2} \phi^\dagger) (\partial^\alpha P_{k_3} \phi).$$
To show that the contribution of this term is $error$, it thus suffices by \eqref{infty-bound} and the triangle inequality to show that
$$ \sum_{k_2, k_3: \max(k_2,k_3) > 10, |k_2 - k_3| \leq 5} \| (\partial_\alpha P_{k_2} \phi^\dagger) (\partial^\alpha P_{k_3} \phi) \|_{L^1_t L^2_x} \lesssim C_0^3 \eps c_0.$$
We use H\"older to split $L^1_t L^2_x$ into two $L^2_t L^4_x$ norms and apply \eqref{four}, \eqref{deriv} to estimate the left-hand side by
$$ \lesssim \sum_{k_2, k_3: \max(k_2,k_3) > 10, |k_2 - k_3| \leq 5} C_0 c_{k_1}
2^{k_2} 2^{-\frac{k_2}{2} - \frac{nk_2}{4}}
C_0 c_{k_2} 2^{k_3} 2^{-\frac{k_3}{2} - \frac{nk_3}{4}}$$
which by \eqref{local} is bounded by
$$ \lesssim C_0^2 c_0^2 \sum_{k_2 > 10} 2^{2k_2(\sigma + 1 - \frac{1}{2} - \frac{n}{4})}.$$
But this is acceptable from our choice of $\sigma$ and the fact that $c_0 \lesssim \eps$.

Now we consider \eqref{dec2}.  By symmetry it suffices to consider the contribution when $k_2 > k_3 + 5, 10$.  In this case we may assume that $|k_1 - k_2| \leq 5$ since the contribution to \eqref{proj} vanishes otherwise.  By the triangle inequality it thus suffices to show that
$$ \sum_{k_2 > 10} \| (P_{k_2 - 5 \leq \cdot \leq k_2 + 5} \phi) (\partial_\alpha P_{k_2} \phi^\dagger) (\partial^\alpha P_{<k_2-5} \phi) \|_{L^1_t L^2_x} \lesssim C_0^3 \eps c_0.$$
We use H\"older, splitting $L^2_t L^4_x$, $L^2_t L^4_x$, $L^\infty_t L^\infty_x$, and use \eqref{four}, \eqref{infty}, \eqref{deriv} and decomposition into projections $P_k$ to estimate the left-hand side by
$$ \sum_{k_2 > 10} (C_0 c_{k_2} 2^{-\frac{k_2}{2} - \frac{nk}{4}})
(C_0 c_{k_2} 2^{k_2} 2^{-\frac{k_2}{2} - \frac{nk}{4}})
(2^{k_2})$$
which is acceptable by the same calculation used to treat \eqref{dec1}.

Now consider \eqref{dec3}.  We may assume that $k_1 < 15$ since the contribution to \eqref{proj} vanishes otherwise.  We can thus simplify \eqref{dec3} as
$$(P_{-10 < \cdot < 15} \phi) (\partial_\alpha P_{<10} \phi^\dagger) (\partial^\alpha P_{<10} \phi).$$
By H\"older it thus suffices to show that
$$\| P_{-10 < \cdot < 15} \phi \|_{L^\infty_t L^2_x} \| \partial_\alpha P_{<10} \phi \|_{L^2_t L^\infty_x} \| \partial^\alpha P_{<10} \phi \|_{L^2_t L^\infty_x} \lesssim C_0^3 \eps c_0.$$
But this is immediate from \eqref{energy}, \eqref{deriv}, and (a trivial modification of) \eqref{2infty-bound}.

The contribution of \eqref{dec4} to \eqref{proj} is always zero, so we turn to \eqref{dec5}.  In light of \eqref{infty-bound} and H\"older it suffices to show that
$$ \| \partial_\alpha P_{-10 < \cdot < 10} \phi \|_{L^2_t L^4_x}^2 \lesssim C_0^3 \eps c_0.$$
But this is immediate from \eqref{four}, \eqref{deriv}, and a breakdown into projections $P_k$.

The terms \eqref{dec6} and \eqref{dec7} are equal. We have thus shown that
$$
P_0(\phi \partial_\alpha \phi^\dagger \partial^\alpha \phi) = 2 
P_0(\tilde \phi \partial_\alpha \tilde \phi^\dagger \partial^\alpha P_{-10 < \cdot < 10} \phi) + error.$$
Since $\psi = P_0 P_{-10 < \cdot < 10} \phi$, it only remains to show the commutator estimate
$$
\| P_0(R \Psi) - R P_0(\Psi)\|_{L^1_t L^2_x} \lesssim C_0^3 \eps c_0
$$
where $R$ is the matrix $R := \tilde \phi \partial_\alpha \tilde \phi^\dagger$ and $\Psi$ is the function $\Psi := \partial^\alpha P_{-10 < \cdot < 10} \phi$.

From \eqref{dual}, \eqref{deriv} and summing over Littlewood-Paley  pieces we have
$$ \| \nabla^2 \tilde \phi \|_{L^2_t L^{n-1}_x} \lesssim C_0 \eps.$$
A similar argument using \eqref{4dual}, \eqref{deriv} gives
$$ \| \nabla \tilde \phi \|_{L^4_t L^{2(n-1)}_x} \lesssim C_0 \eps.$$
Combining these together with \eqref{infty-bound} we obtain
\be{stew}
\| \nabla R \|_{L^2_t L^{n-1}_x} \lesssim C^2_0 \eps
\end{equation}
while from \eqref{endpoint}, \eqref{deriv} we have
$$ \| \Psi \|_{L^2_t L^{2(n-1)/(n-3)}_x} \lesssim C_0 c_0.$$
Thus to finish the proof of this Proposition it suffices to use the standard commutator estimate (with $p = n-1$, $q = 2(n-1)/(n-3)$, and $r=2$)

\begin{lemma}\label{cook-lemma}  We have
\be{cook}
\| P_0(fg) - f P_0(g) \|_r \lesssim \| \nabla f\|_{p} \|g\|_{q}
\end{equation}
for all smooth functions $f$, $g$ on $\R^n$ and all $1 \leq p, q, r \leq \infty$ such that $1/p + 1/q = 1/r$.  
\end{lemma}

For a previous application of this type of lemma to wave maps, see \cite{keel:wavemap}.

\begin{proof}
We begin with the identity
$$ P_0(fg)(x) - f P_0(g)(x) = \int \check m(y) (f(x-y) - f(x)) g(x-y)\ dy$$
and use the Fundamental theorem of Calculus to rewrite this as
$$- \int_0^1 \int \check m(y) y \cdot \nabla f(x-ty) g(x-y)\ dy\ dt.$$
Since $\check m(y) y$ is integrable, the claim then follows from Minkowski
and H\"older.
\end{proof}

\end{proof}

We have thus transformed the non-linear equation \eqref{wavemap-eq} into the linear equation \eqref{proj-effect} (modulo a forcing term which can be dealt with).  This procedure appears to be quite general, and does not rely on the specific form of \eqref{wavemap-eq}.  In principle, the linearity conveys several advantages, for instance we can now apply Duhamel's principle to remove the inhomogeneity, or the principle of superposition to decompose $\psi$ into wave packets or similar objects.  However, we shall not exploit linearity in any significant way; for us, the main advantage of \eqref{proj-effect} is that there is only one moderately high-frequency term in the non-linearity, with the other two factors being very smooth.  Also, the linearity allows us to separate the role of the small quantities $\eps$ and $c_0$; $\psi$ will always be estimated with a bound of $c_0$, whereas $\tilde \phi$ will be estimated with a factor of $\eps$ if it is high-frequency or is accompanied with a derivative, or by a constant otherwise.

The main term in the right-hand side of \eqref{proj-effect} cannot be placed in $L^1_t L^2_x$ by Strichartz estimates; in fact the failure of the Strichartz estimates here is quite dramatic (more than a logarithm).  Roughly speaking, this is because the derivatives in this term could fall on high-frequency components of $\phi$ rather than low frequency ones.  However, after the renormalization we will be able to move all the derivatives onto low frequency terms, allowing Strichartz estimates to successfully place the non-linearity into $L^1_t L^2_x$.

\section{The tangent space of the sphere}\label{tangent-sec}

We now develop the analogue of \eqref{nabla}.

\begin{proposition}\label{sphere}  We have the estimate
\be{sph1}
\|\tilde \phi^\dagger \partial^\alpha \psi \|_{L^2_t L^2_x} \lesssim C_0^2 \eps c_0
\end{equation}
\end{proposition}

\begin{proof}
We first observe that \eqref{sph1} is easy to show if the derivative is moved to the low frequency term:
\be{oomp}
\|(\partial^\alpha \tilde \phi^\dagger) \psi \|_{L^2_t L^2_x} \lesssim C_0^2 \eps c_0.
\end{equation}
Indeed, we simply place $\psi$ in $L^\infty_t L^2_x$ using \eqref{energy}, \eqref{deriv} and $\partial^\alpha \tilde \phi$ in $L^2_t L^\infty_x$ using \eqref{2infty}, \eqref{deriv}.  From this and the product rule it thus suffices to show that
$$
\|\partial^\alpha (\tilde \phi^\dagger \psi) \|_{L^2_t L^2_x} \lesssim C_0^2 \eps c_0.$$
Consider the expression $P_0(\phi^\dagger \phi)$.  Heuristically, this quantity is approximately $2 \tilde \phi^\dagger \psi$, by the same type of reasoning used to obtain the linearization \eqref{linearized}.  On the other hand, since $\phi$ lies on the sphere, $P_0(\phi^\dagger \phi) = P_0(1) = 0$.  Accordingly, we shall rewrite the above estimate as
\be{choo}
\|\partial^\alpha (P_0(\phi^\dagger \phi) - 2\tilde \phi^\dagger \psi) \|_{L^2_t L^2_x} \lesssim C_0^2 \eps c_0.
\end{equation}
We now split $\phi^\dagger \phi$ as
$$
(P_{>-10} \phi^\dagger) (P_{>-10} \phi) + \tilde \phi^\dagger (P_{-10 < \cdot < 10} \phi) + (P_{-10 < \cdot < 10} \phi^\dagger) \tilde \phi
$$
plus other terms which vanish when $P_0$ is applied.

The first term of the above expression can be refined to
$$ \sum_{k_1, k_2 > -10: |k_1 - k_2| < 20} (P_{k_1} \phi^\dagger) (P_{k_2} \phi)$$
since all the other components vanish after applying $P_0$.  The contribution of this term to \eqref{choo} can thus be estimated using the product rule and H\"older by
$$ \lesssim \sum_{k_1,k_2 > -10} \sum_{k_2: |k_1 - k_2| < 20} \| \partial^\alpha P_{k_1} \phi \|_{L^\infty_t L^2_x} \| P_{k_2} \phi \|_{L^2_t L^\infty_x}
+ \| P_{k_1} \phi \|_{L^2_t L^\infty_x} \| \partial^\alpha P_{k_2} \phi \|_{L^\infty_t L^2_x}.$$
Applying \eqref{2infty}, \eqref{energy}, \eqref{deriv}, this can be bounded by
$$ \lesssim \sum_{k_1 > -10} 2^{k_1} C_0 c_{k_1} 2^{-nk_1/2} C_0 c_{k_1} 2^{-k_1/2}$$
which is acceptable by \eqref{local} and our choice of $\delta$.

The other two terms are equal to each other.  It thus remains to show the commutator estimate
$$ \| \partial^\alpha( P_0(\tilde \phi^\dagger \Psi) - \tilde \phi^\dagger P_0(\Psi) ) \|_{L^2_t L^2_x} \lesssim C_0^2 \eps c_0$$
where $\Psi := P_{-10 < \cdot < 10} \phi$.  Since the expression inside the $\partial^\alpha$ has Fourier support on $|\xi| \lesssim 1$, we may discard the derivative $\partial^\alpha$.

By Lemma \ref{cook-lemma} we may estimate the left-hand side of this by
$$ \| \nabla \tilde \phi\|_{L^2_t L^\infty_x} \| \Psi \|_{L^\infty_t L^2_x},$$
and this is acceptable by the argument used to treat \eqref{oomp} (as the commutator estimate has effectively moved the derivative from the high-frequency term to the low-frequency term; cf. the ``I-method'' trick in \cite{keel:wavemap}, \cite{keel:mkg}, \cite{iteam-1}).
\end{proof}

As a particular corollary of \eqref{sph1}, \eqref{2infty-bound}, and H\"older we have
$$
\partial_\alpha \tilde \phi \tilde \phi^\dagger \partial^\alpha \psi = error$$
and so from \eqref{proj-effect} we have the analogue 
\be{cancel-accurate}
\Box \psi = 2 A_\alpha \partial^\alpha \psi + error
\end{equation}
of \eqref{cancel}, where $A_\alpha$ is the anti-symmetric $m \times m$ matrix
$$ A_\alpha := \partial_\alpha \tilde \phi \tilde \phi^\dagger - \tilde \phi \partial_\alpha \tilde \phi^\dagger.$$

The equation \eqref{cancel-accurate} can be derived in an alternate way, which we now sketch.  By differentiating \eqref{wavemap-eq} we have
$$ \Box \nabla_{x,t} \phi = - 2 \phi \partial_\alpha \phi^\dagger \partial^\alpha \nabla_{x,t} \phi - \nabla \phi_{x,t} \partial_\alpha \phi^\dagger \partial^\alpha \phi.$$
Since $\phi$ lies on the sphere, we have $\phi^\dagger \nabla_{x,t} \phi = 0$, so we can rewrite this as
$$ \Box \nabla_{x,t} \phi = - 2 A_{\alpha,\infty} \partial^\alpha \nabla_{x,t} \phi - \nabla_{x,t} \phi \partial^\alpha \phi^\dagger \partial^\alpha \phi$$
where
$$ A_{\alpha,\infty} := \partial_\alpha \phi \phi^\dagger - \phi \partial_\alpha \phi^\dagger.$$
If one then applies a Littlewood-Paley projection $P_k$ to this equation then one can obtain (the derivative of) \eqref{cancel-accurate} by the same type of calculations as in the previous section.  Note that the second term in the above equation is cubic in $\nabla_{x,t} \phi$ and so can be treated by Strichartz estimates since $n \geq 3$.  This derivation of \eqref{cancel-accurate} works well in high dimensions but is difficult to use in low dimensions, especially $n=2$, because one only has $n/2$ degrees of regularity, and so we cannot afford to differentiate the equation as above.

\section{Approximate parallel transport}\label{transport-sec}

We now construct a matrix field $U$ which is approximately orthogonal and which will renormalize \eqref{cancel-accurate} into a much better form, namely $\Box w = error$.

We shall use the scheme described in the introduction.  More precisely, we let $M$ be a large integer (depending on $T$!) to be chosen later  and define the real $m \times m$-valued matrix field $U$ by
$$ U := I + \sum_{-M < k \leq -10} U_k$$ 
where $I$ is the identity matrix and the $U_k$ are defined inductively by
\be{ack-accurate}
U_k := ((P_k \phi) (P_{<k} \phi^\dagger) - (P_{<k} \phi) (P_k \phi^\dagger)) U_{<k}
\end{equation}
and
$$ U_{<k} := I + \sum_{-M < k' < k} U_{k'}.$$

An easy inductive argument shows that $U_{<k}$ has Fourier support on the region $\{ |\xi| \leq 2^{k+5} \}$, and thus that $U$ has Fourier support on the region $\{ |\xi| \leq 2^{-5} \}$.  More generally, one expects $U_{<k}$ to have essentially the same estimates as $P_{<k} \phi$, using the heuristic that the factors $P_{<k} \phi$ and $U_{<k}$ in \eqref{ack-accurate} are bounded and therefore do not significantly affect the estimates. 

We now quantify the precise estimates on $U$ which we shall need. 

\begin{proposition}\label{u-est}
Assume that $\eps$ is sufficiently small depending on $C_0$, and $M$ is sufficiently large depending on $T$, $C_0$, $\eps$.  Then we have the almost orthogonality property
\be{u-invert}
\| U^\dagger U - I \|_{L^\infty_t L^\infty_x}, 
\| \partial_t(U^\dagger U - I) \|_{L^\infty_t L^\infty_x}
 \lesssim C_0^2 \eps.
\end{equation}
In particular, if $\eps$ is sufficiently small depending on $C_0$, then $U$ is invertible, and
\be{u-infty}
\| U \|_{L^\infty_t L^\infty_x}, \| U^{-1} \|_{L^\infty_t L^\infty_x} \lesssim 1.
\end{equation}
Also, we have the approximate parallel transport property (cf. \eqref{ook})
\be{1infty}
\| \partial_\alpha U - A_\alpha U \|_{L^1_t L^\infty_x} \lesssim C_0^2 \eps
\end{equation}
as well as the additional bounds (needed to control error terms)
\be{u-energy}
\| \partial_\alpha U \|_{L^\infty_t L^\infty_x} \lesssim C_0^2 \eps
\end{equation}
\be{u2}
\| \partial_\alpha U \|_{L^2_t L^\infty_x} \lesssim C_0^2 \eps
\end{equation}
\be{u-dd}
\| \Box U \|_{L^2_t L^{n-1}_x} \lesssim C_0^2 \eps
\end{equation}
for all $\alpha$.
\end{proposition}

The power of $\eps$ is not sharp in many of these estimates, but that is irrelevant for our purposes.  Interestingly, the above bounds on $U$ do not seem to easily extend to any useful continuity estimates on the map $\phi \mapsto U$.  In particular, small $\dot H^{n/2}$ perturbations in $\phi$ can lead to large fluctuations in $U$ in $L^\infty$.  This phenomenon is the major obstacle to obtaining a critical Sobolev regularity well-posedness theory from our arguments, and is also a problem in obtaining scattering even for classical wave maps.

\begin{proof}
As noted in the introduction, we have the identity
$$ 
U_k^\dagger U_{<k}+  U_{<k}^\dagger U_k = 0$$
whence 
\be{uident}
U^\dagger_{<K} U_{<K} - I = \sum_{-M < k < K} U_k^\dagger U_k
\end{equation}
for all $-M \leq K \leq -9$.

We now show inductively that
\be{induct}
\| U_{<K} \|_{L^\infty_t L^\infty_x} \leq 2
\end{equation}
for all $-M \leq K \leq -9$.
This is clearly true for $K=-M$.  Now suppose that $K > -M$ and the claim has been proven for all smaller $K$.  Then from \eqref{ack-accurate} and H\"older we have
$$ \| U_k \|_{L^\infty_t L^\infty_x} \lesssim \| P_k \phi \|_{L^\infty_t L^\infty_x} \| P_{<k} \phi \|_{L^\infty_t L^\infty_x}$$
for all $-M < k < K$.  By \eqref{infty}, \eqref{deriv}, and \eqref{infty-bound} we thus have 
$$ \| U_k \|_{L^\infty_t L^\infty_x} \lesssim C_0 c_k,$$
and the induction \eqref{induct} can thus be closed by \eqref{uident} and \eqref{l2-ass}, if $\eps$ is sufficiently small.

From the above analysis we see that the first part of \eqref{u-invert} obtains, as does the first part of \eqref{u-infty}.  In particular $U^\dagger U$, and thus $U$, are invertible, and this gives the second part of \eqref{u-infty}.

We now show \eqref{u-energy}, \eqref{u2}; the second part of \eqref{u-invert} will then follow from \eqref{u-infty}, \eqref{u-energy}, and the product rule.

We shall again use induction, showing that
\be{choco}
\| \partial_\alpha U_{<K} \|_{L^\infty_t L^\infty_x} \leq C_1 2^K C_0^2 c_K
\end{equation}
and
\be{chico}
\| \partial_\alpha U_{<K} \|_{L^2_t L^\infty_x} \leq C_1 2^{K/2} C_0^2 c_K
\end{equation}
for all $-M \leq K \leq -9$ and some sufficiently large absolute constant $C_1 \gg 1$.

The claim is trivial when $K=-M$.  Now suppose that $K > M$ and the claim has been proven for all smaller $K$.  By differentiating \eqref{ack-accurate} and using H\"older, we obtain
\begin{align*}
\| \partial_\alpha U_{K-1} \|_{L^\infty_t L^\infty_x} \lesssim
&\| \partial_\alpha P_{K-1} \phi \|_{L^\infty_t L^\infty_x} \| P_{<K-1} \phi \|_{L^\infty_t L^\infty_x} \| U_{<K-1} \|_{L^\infty_t L^\infty_x}\\
&+
\| P_{K-1} \phi \|_{L^\infty_t L^\infty_x} \| \partial_\alpha P_{<K-1} \phi \|_{L^\infty_t L^\infty_x} \| U_{<K-1} \|_{L^\infty_t L^\infty_x}\\
&+
\| P_{K-1} \phi \|_{L^\infty_t L^\infty_x} \| P_{<K-1} \phi \|_{L^\infty_t L^\infty_x} \| \partial_\alpha U_{<K-1} \|_{L^\infty_t L^\infty_x}
\end{align*}
and
\begin{align*}
\| \partial_\alpha U_{K-1} \|_{L^2_t L^\infty_x} \lesssim
&\| \partial_\alpha P_{K-1} \phi \|_{L^2_t L^\infty_x} \| P_{<K-1} \phi \|_{L^\infty_t L^\infty_x} \| U_{<K-1} \|_{L^\infty_t L^\infty_x}\\
&+
\| P_{K-1} \phi \|_{L^\infty_t L^\infty_x} \| \partial_\alpha P_{<K-1} \phi \|_{L^2_t L^\infty_x} \| U_{<K-1} \|_{L^\infty_t L^\infty_x}\\
&+
\| P_{K-1} \phi \|_{L^\infty_t L^\infty_x} \| P_{<K-1} \phi \|_{L^\infty_t L^\infty_x} \| \partial_\alpha U_{<K-1} \|_{L^2_t L^\infty_x}.
\end{align*}
By applying \eqref{infty}, \eqref{2infty}, \eqref{deriv}, \eqref{u-infty}, and the induction hypothesis we thus see that
$$ 
\| \partial_\alpha U_{K-1} \|_{L^\infty_t L^\infty_x} \lesssim
2^K C_0^2 c_K (1 + C_1\eps)$$
and
$$ 
\| \partial_\alpha U_{K-1} \|_{L^2_t L^\infty_x} \lesssim
2^{K/2} C_0^2 c_K (1 + C_1\eps).$$
If $\eps$ is sufficiently small and $C_1$ is sufficiently large depending on $C_0$  then one can close the induction hypothesis.  This gives \eqref{u-energy}, \eqref{u2}.

Next, we prove \eqref{1infty}.  We can write $A_\alpha$ as $A_{\alpha, \leq -10}$ where
$$ A_{\alpha, \leq k} := A_{\alpha, <k+1} := (\partial_\alpha P_{\leq k} \phi) P_{\leq k} \phi^\dagger - (P_{\leq k} \phi) (\partial_\alpha P_{\leq k} \phi)^\dagger.$$
We thus have the telescoping identity
$$ \partial_\alpha U - A_\alpha U = [\sum_{-M < k \leq -10} \partial_\alpha U_k - (A_{\alpha,\leq k} U_{\leq k} - A_{\alpha,<k} U_{<k})] - A_{\alpha,\leq -M}.$$
To estimate $A_{\alpha,\leq -M}$, we use \eqref{infty}, \eqref{deriv} to obtain
$$ \| A_{\alpha,\leq -M} \|_{L^\infty_t L^\infty_x} \lesssim 2^{-M}.$$
This term is thus acceptable by H\"older in time if $M$ is sufficiently large depending on $T$, $C_0$, $\eps$.

By \eqref{l2-ass} and the triangle inequality it thus suffices to show that
$$
\| \partial_\alpha U_k - (A_{\alpha,\leq k} U_{\leq k} - A_{\alpha,<k} U_{<k})\|_{L^1_t L^\infty_x} \lesssim C_0^2 c_k^2$$
for all $-M < k \leq -10$.

We expand out
\begin{align}
\partial_\alpha U_k &= 
((\partial_\alpha P_k \phi) (P_{<k} \phi^\dagger) - (P_{<k} \phi) (\partial_\alpha P_k \phi^\dagger)) U_{<k}\label{d1}\\
&+
((P_k \phi) (\partial_\alpha P_{<k} \phi^\dagger) - (\partial_\alpha P_{<k} \phi) (P_k \phi^\dagger)) U_{<k}\nonumber\\
&+
((P_k \phi) (P_{<k} \phi^\dagger) - (P_{<k} \phi) (P_k \phi^\dagger)) \partial_\alpha U_{<k}\nonumber
\end{align}
and
\begin{align}
A_{\alpha,\leq k} U_{\leq k} - A_{\alpha,<k} U_{<k} &=
((\partial_\alpha P_k \phi) (P_{<k} \phi^\dagger) - (P_{<k} \phi) (\partial_\alpha P_k \phi^\dagger)) U_{<k}\label{d2}\\
&+
((\partial_\alpha P_{\leq k} \phi) (P_k \phi^\dagger) - (P_k \phi) (\partial_\alpha P_{\leq k} \phi^\dagger)) U_{<k}\nonumber\\
&+
((\partial_\alpha P_{\leq k} \phi) (P_{\leq k} \phi^\dagger) - (P_{\leq k} \phi) (\partial_\alpha P_{\leq k} \phi^\dagger)) U_{k}\nonumber.
\end{align}
In both expressions, the dangerous terms \eqref{d1}, \eqref{d2} occur when the derivative falls on a high frequency term $P_k \phi$ instead of a low frequency term such as $P_{<k} \phi$, $P_{\leq k} \phi$, $U_{<k}$.  (Indeed \eqref{d2} is the only reason why $A_\alpha$ fails to be in $L^1_t L^\infty_x$, and is the only reason why we need a renormalization by $U$ in the first place).  Fortunately, we have chosen $U$ so that the dangerous terms \eqref{d1}, \eqref{d2} cancel each other.  From the triangle inequality it thus suffices to show the bounds
\be{chunk4}
\begin{split}
\| |P_k \phi| |\partial_\alpha P_{<k} \phi| |U_{<k}| \|_{L^1_t L^\infty_x},
\| |P_k \phi| |P_{<k} \phi| |\partial_\alpha U_{<k}| \|_{L^1_t L^\infty_x},&\\
\| |\partial_\alpha P_{\leq k} \phi| |P_k \phi| |U_{<k}| \|_{L^1_t L^\infty_x},
\| |\partial_\alpha P_{\leq k} \phi| |P_{\leq k} \phi| |U_k| \|_{L^1_t L^\infty_x}&
\lesssim C_0^2 c_k^2.
\end{split}
\end{equation}
The first term in \eqref{chunk4} is acceptable by \eqref{2infty}, \eqref{deriv} for the first two factors (dyadically decomposing the latter factor and using \eqref{local}) and \eqref{u-infty} for the last.  The second term is acceptable by \eqref{2infty}, \eqref{deriv} for the first factor, \eqref{infty-bound} for the second, and \eqref{chico} for the last.  The third term is treatable by the same argument as the first term.  Finally, the fourth term is acceptable by \eqref{2infty}, \eqref{deriv} for the first term, \eqref{infty-bound} for the second term, and the estimate
$$
\| U_k \|_{L^2_t L^\infty_x} \lesssim 2^{-k/2} C_0 c_k,$$
which can be proven from \eqref{ack-accurate} and estimating $P_k \phi$ using \eqref{2infty}, \eqref{deriv}, $P_{<k} \phi$ using \eqref{infty-bound}, and $U_{<k}$ using \eqref{u-infty}.

The only remaining estimate to prove is \eqref{u-dd}.  In principle this is the same type of estimate as \eqref{stew}, but there is a minor complication arising from the double time derivative in $\Box$, which are not directly treatable by the $S_k$ norms.  To get around this we will have to use the equation \eqref{wavemap-eq}.

More precisely, we shall need 

\begin{lemma}\label{box}  For all $k$, we have
\be{box-est}
\| \Box P_k \phi \|_{L^2_t L^{n-1}_x} \lesssim 2^{k(2 - \frac{1}{2} - \frac{n}{n-1})} C_0^3 c_k.
\end{equation}
\end{lemma}

\begin{proof}
Morally speaking this estimate obtains from \eqref{dual}, \eqref{deriv} if we treat time derivatives like spatial ones.  We could have modified our Littlewood-Paley  operators to project in time as well as space in order to make this heuristic rigorous, but this creates other difficulties having to do with time localization which we wished to avoid.

We shall show \eqref{box-est} for $k=0$ to simplify the exposition; the reader may verify that the argument below is scale invariant and thus extends to all $k$.

Applying \eqref{wavemap-eq}, we see it suffices to show that
$$ \| P_0( \phi \partial_\alpha \phi^\dagger \partial_\alpha \phi ) \|_{L^2_t L^{n-1}_x} \lesssim C_0^3 c_0.$$

Let us first consider the contribution of
$$ \phi \partial_\alpha P_{>5} \phi^\dagger \partial_\alpha P_{>10} \phi.$$
In this case it suffices by Bernstein's inequality (or Young's inequality) to obtain $L^2_t L^2_x$ estimates.  From \eqref{deriv} and the definition of the $S_k$ norm we have
$$ \| \partial_\alpha P_k \phi \|_{L^4_t L^4_x} \lesssim C_0 c_k 2^{-(n-3)k/4}.$$
From \eqref{local} and the assumptions on $\sigma$ we thus have
$$ \| \partial_\alpha P_{>5} \phi \|_{L^4_t L^4_x}, \| \partial_\alpha P_{>10} \phi \|_{L^4_t L^4_x} \lesssim C_0 c_0$$
and this contribution is thus acceptable by \eqref{infty-bound}.

Now consider the contribution of
$$ \phi \partial_\alpha P_{\leq 5} \phi^\dagger \partial_\alpha P_{>10} \phi.$$
In this case we can replace the first factor by $P_{>5} \phi$, since the error in doing so vanishes after applying $P_0$.  We now modify the above argument, the only difference being that we now place the first term in $L^4_t L^4_x$ and the second in $L^\infty_t L^\infty_x$.  In fact, the summation is much better because the derivative is now on the low frequency term.

A similar argument deals with
$$ \phi \partial_\alpha P_{> 5} \phi^\dagger \partial_\alpha P_{\leq 10} \phi$$
and so we are left with
$$ \phi \partial_\alpha P_{\leq 5} \phi^\dagger \partial_\alpha P_{\leq 10} \phi.$$
We can split this into
$$ P_{-5 \leq \cdot \leq 15} \phi \partial_\alpha P_{\leq 5} \phi^\dagger \partial_\alpha P_{\leq 10} \phi,$$
$$ P_{<-5} \phi \partial_\alpha P_{-5 \leq \cdot \leq 5} \phi^\dagger \partial_\alpha P_{\leq 10} \phi,$$
and
$$ P_{<-5} \phi \partial_\alpha P_{<-5} \phi^\dagger \partial_\alpha P_{-5 \leq \cdot \leq 10} \phi,$$
plus some other terms which vanish when $P_0$ is applied.

For each of these three terms we place the high frequency ($>-5$) factor in $L^2_t L^{n-1}_x$ and the other two factors in $L^\infty_t L^\infty_x$.  Regardless of the position of the derivatives, the high frequency factor has a norm of $O(C_0 c_0)$ by \eqref{dual}, \eqref{deriv}.  Of the other two factors, both are bounded by $O(1)$ by \eqref{infty-bound}, and at least one contains a derivative and therefore has a norm of $O(C_0 \eps)$ by \eqref{infty}, \eqref{deriv}.  The claim then follows (if $\eps$ is sufficiently small depending on $C_0$).
\end{proof}

We now prove \eqref{u-dd}.  As before, we shall use induction and in fact prove 
$$ \| \Box U_{<K} \|_{L^2_t L^{n-1}_x} \leq C_2 2^{K(2 - \frac{1}{2} - \frac{n}{n-1})} C_0^2 \eps$$
for all $-M \leq K \leq -10$ and some large $C_2 \gg 1$.

The claim is trivial for $K=-M$.  Now suppose that $K > M$ and the claim has been proven for all smaller $K$.  We apply $\Box$ to \eqref{ack-accurate} and take absolute values (ignoring any possibility of cancellation) to obtain
\begin{align*}
|\Box U_{K-1}| \lesssim &|P_{K-1} \phi| |P_{<K-1} \phi| |\Box U_{<K-1}|
+ |P_{K-1} \phi| |\Box P_{<K-1} \phi| |U_{<K-1}|\\
&+ |\Box P_{K-1} \phi| |P_{<K-1} \phi| |U_{<K-1}|
+ |\nabla_{x,t} P_{K-1} \phi| |\nabla_{x,t} P_{<K-1} \phi| |U_{<K-1}|\\
&+ |\nabla_{x,t} P_{K-1} \phi| |P_{<K-1} \phi| |\nabla_{x,t} U_{<K-1}|
+ |P_{K-1} \phi| |\nabla_{x,t} P_{<K-1} \phi| |\nabla_{x,t} U_{<K-1}|.
\end{align*}
We will show that all six terms on the right-hand side have an $L^2_t L^{n-1}_x$ norm of 
$$ \lesssim 2^{K(2 - \frac{1}{2} - \frac{n}{n-1})} C_0^4 \eps (1 + C_2 \eps),$$
so that we can close the induction if $\eps$ is sufficiently small and $C_2$ sufficiently large.

For the first term we use \eqref{infty}, \eqref{deriv} for the first factor, \eqref{infty-bound} for the second factor, and the induction hypothesis for the third factor.  For the second and third terms we use Lemma \ref{box} for the $\Box \phi$ term, and place the other two terms in $L^\infty_t L^\infty_x$ using \eqref{infty-bound}, \eqref{u-infty}.  For the remaining three terms we place the first factor in $L^2_t L^{n-1}_x$ using \eqref{2infty}, \eqref{deriv} and the other two factors in $L^\infty_t L^\infty_x$ using \eqref{infty}, \eqref{deriv} for the second factor and \eqref{choco} for the third.
\end{proof}

\section{Wrapping up}\label{conclusion-sec}

Armed with Proposition \ref{u-est} we can now conclude the proof of \eqref{psi-bound} and thus of Propositions \ref{reduced} and \ref{main-prop}.
Since $U$ is invertible, we may write $\psi = Uw$ for some $w$, which is smooth by our assumptions.  By \eqref{u-infty}, \eqref{u-energy} and the Leibnitz rule in time we have
$$ \| \psi \|_{S_0} \lesssim \| w \|_{S_0}$$
and so it suffices to show that
\be{w-targ}
\|w \|_{S_0} \ll C_0 c_0.
\end{equation}

We expand \eqref{cancel-accurate} using the Leibnitz rule as
$$U \Box w + 2 \partial_\alpha U \partial^\alpha w + (\Box U) w = 2 A_\alpha U \partial^\alpha w + 2 A_\alpha (\partial^\alpha U) w + error.$$
By \eqref{u-infty} we see that $U^{-1} error = error$, thus we can rewrite the previous as
\be{chunk}
\Box w = -2 U^{-1} (\partial_\alpha U - A_\alpha U) \partial^\alpha w + 2 U^{-1} A_\alpha (\partial^\alpha U) U^{-1} \psi - U^{-1} (\Box U) U^{-1} \psi + error.
\end{equation}
We now show that all terms on the right-hand side are of the form $error$. 

To control the first term, it suffices by \eqref{u-infty}, \eqref{1infty} (if $\eps$ is sufficiently small depending on $C_0$) to show that
$$ \| \partial^\alpha w \|_{L^\infty_t L^2_x} \lesssim C_0 c_0.$$

Since $\psi = Uw$ and $\partial^\alpha \psi = U \partial^\alpha w + (\partial^\alpha U) w$, we have
$$ \partial^\alpha w = U^{-1} \partial^\alpha \psi + U^{-1} (\partial^\alpha U) U^{-1}\psi.$$
The claim then follows from \eqref{u-infty}, \eqref{energy}, \eqref{deriv} (for the first term) and \eqref{u-infty}, \eqref{u-energy}, \eqref{energy}, \eqref{deriv} (for the second term).

Now consider the second term in \eqref{chunk}.  By \eqref{infty-bound}, \eqref{2infty-bound} we have
$$ \| A_\alpha \|_{L^2_t L^\infty_x} \lesssim C_0 \eps$$
and from \eqref{2infty}, \eqref{energy} we have
$$ \| \psi \|_{L^\infty_t L^2_x} \lesssim C_0 c_0.$$
The claim then follows from \eqref{u-infty} and \eqref{u2}.

Now consider the third term in \eqref{chunk}.  By \eqref{endpoint}, \eqref{deriv} we have
$$ \| \psi \|_{L^2_t L^{2(n-1)/(n-3)}} \lesssim C_0 c_0.$$
The claim then follows from \eqref{u-dd} and \eqref{u-infty}.

We have thus shown that $\Box w = error$, or in other words that
$$ \| \Box w\|_{L^1_t L^2_x} \lesssim C_0^3 \eps c_0$$
Also, from \eqref{u-infty} and the assumption on $\psi[0]$ we have
$$ \| w[0]\|_{L^2_x} \lesssim c_0.$$
At this point one should be able to obtain \eqref{w-targ} from Theorem \ref{strich}, however $w$ is not quite supported on the frequency annulus $|\xi| \sim 1$ (we have Fourier support control on $U$, but not on $U^{-1}$).  However, we can apply $P_{-10 < \cdot < 10}$ to the above estimates and use Theorem \ref{strich} to conclude that
\be{nearly}
\| P_{-10 < \cdot < 10} w \|_{S_0} \ll C_0 c_0
\end{equation}
(providing that $C_0$ is sufficiently large depending only on $n$, $m$, $\sigma$, and $\eps$ is sufficiently small depending on $C_0$).  

To pass from \eqref{nearly} to \eqref{w-targ} we begin with the identity
$$ w = U^\dagger \psi - (U^\dagger U - I) w.$$
From the Fourier support of $U$ and $\psi$ we see that $U^\dagger \psi$ has Fourier support in the annulus $2^{-5} \leq |\xi| \leq 2^5$, which implies that
$$ (1 - P_{-10 < \cdot < 10}) w = -(1 - P_{-10 < \cdot < 10}) (U^\dagger U - I) w.$$
Taking $S_0$ norms of both sides we see that
$$ \| (1 - P_{-10 < \cdot < 10}) w \|_{S_0} \lesssim \| (U^\dagger U - I) w \|_{S_0}.$$
From \eqref{u-invert}, and H\"older we thus obtain
$$ \|w - P_{-10 < \cdot < 10} w \|_{S_0}  \lesssim C_0^2 \eps^2 \|w\|_{S_0}$$
and \eqref{w-targ} then follows from \eqref{nearly} if $C_0$ is sufficiently large and $\eps$ is sufficiently small.  This concludes the proof of Proposition \ref{main-prop}, and thus of Theorem \ref{main2}.

\endprf


\begin{thebibliography}{10}

\bibitem{changwangyang}
S.Y.A. Chang, L. Wang, P. Yang, \emph{Regularity of Harmonic Maps},
Comm. Pure. Appl. Math. \textbf{52} (1999), 1099--1111.

\bibitem{christ.spherical.wave}
D. Christodoulou, A. Tahvildar-Zadeh, \emph{On the regularity of 
spherically symmetric wave maps}, Comm. Pure Appl. Math, \textbf{46} (1993),
1041--1091.

\bibitem{fms:compact}
A. Freire, S. M\"uller, M. Struwe, \emph{Weak compactness of wave maps and harmonic maps}, Ann. Inst. H. Poincare Anal. Non Lineaire \textbf{15} (1998), no. 6, 725--754.

\bibitem{iteam-1}
J. Colliander, M. Keel, G. Staffilani, H. Takaoka, T. Tao, \emph{Global well-posedness below the energy norm for 2D NLS}, preprint.

\bibitem{ginebre}
J. Ginibre, G. Velo, \emph{The Cauchy problem for the $O(N)$, $CP(N-1)$,
and $GC(N,P)$ models},
Ann. Physics, \textbf{142} (1982), 
393--415.

\bibitem{grillakis.equivariant}
M. Grillakis, \emph{Classical solutions for the equivariant wave map in 
$1+2$ dimensions}, to appear in Indiana Univ. Math. J.

\bibitem{gu}
C. Gu, \emph{On the Cauchy problem for harmonic maps defined on two-dimensional
Minkowski space}, 
Comm. Pure Appl. Math., \textbf{33},(1980), 
727--737.


\bibitem{helein}
F. Helein, \emph{Regularite des applications faiblement harmoniques entre une sur
face et une varitee Riemannienne}, C.R. Acad. Sci. Paris Ser. I Math.,
\textbf{312} (1991), 591-596.

\bibitem{tao:keel}
M. Keel, T. Tao, \emph{Endpoint Strichartz Estimates}, Amer. Math. J. 120 (1998), 955--980.

\bibitem{keel:wavemap}
M. Keel, T. Tao, \emph{Local and global well-posedness
of wave maps on $\R^{1+1}$ for rough data}, IMRN \textbf{21} (1998),
1117--1156.

\bibitem{keel:mkg}
M. Keel, T. Tao, \emph{Global existence for the Maxwell-Klein-Gordon equation below the energy norm}, in preparation.

\bibitem{kman.barrett}
S. Klainerman, \emph{On the regularity of classical field theories in Minkowski
space-time $\R^{3+1}$}, Prog. in Nonlin. Diff. Eq. and their Applic., \textbf{
29},
(1997), Birkh\"auser, 113--150.

\bibitem{kman.mach.smoothing}
S. Klainerman, M. Machedon, \emph{Smoothing estimates for null forms
and applications}, Duke Math. J., \textbf{81} (1995), 99--133.

\bibitem{kl-mac:null3}
S. Klainerman, M. Machedon, \emph{On the optimal local regularity for gauge field theories,} Diff. and Integral Eq. \textbf{10} (1997), 1019--1030.

\bibitem{kman.selberg}
S. Klainerman, S. Selberg, \emph{ Remark on the optimal regularity for
equations of wave maps type}, C.P.D.E., \textbf{22} (1997), 901--918. 

\bibitem{kman.selberg:survey}
S. Klainerman, S. Selberg, \emph{Bilinear estimates and applications to nonlinear wave equations}, preprint. 

\bibitem{lady}
O.A. Ladyzhenskaya, V.I. Shubov, \emph{Unique solvability of the Cauchy
problem for the equations of the two dimensional chiral fields, taking values
in complete Riemann manifolds}, J. Soviet Math.,
\textbf{25} (1984), 855--864. (English Trans. of 1981 Article.)

\bibitem{nak:wavemap}
K. Nakanishi, \emph{Local well-posedness and Illposedness in the critical Besov spaces for semilinear wave equations with quadratic forms}, Funk. Ekvac. \textbf{42} (1999), 261-279.

\bibitem{selberg:thesis}
S. Selberg, \emph{Multilinear space-time estimates and applications to local existence theory for non-linear wave equations}, Princeton University Thesis.

\bibitem{selberg:wavemaps}
S. Selberg, \emph{Wave maps and bilinear spacetime estimates}, preprint.

\bibitem{shatah}
J. Shatah, \emph{Weak solutions and development of singularities
of the $SU(2)$ $\sigma$-model.}
Comm. Pure Appl. Math., \textbf{41} (1988),
459--469.

\bibitem{shatah.zurich}
J. Shatah, \emph{The Cauchy problem for harmonic maps on Minkowski space},
in Proceed. Inter. Congress of Math. 1994, Birkh\"auser, 1126--1132.

\bibitem{shatah-struwe}
J. Shatah, M. Struwe, \emph{Geometric Wave Equations}, Courant Lecture Notes in Mathematics \textbf{2} (1998)

\bibitem{shatah.shadi.blow} 
J. Shatah, A. Tavildar-Zadeh, \emph{On the Cauchy problem for equivariant
wave maps}, Comm. Pure Appl. Math., \textbf{47} (1994), 719 - 753.

\bibitem{sideris.harmonic}
T. Sideris, \emph{Global existence of harmonic maps in Minkowski space},
Comm. Pure Appl. Math., \textbf{42} (1989),1--13.
 
\bibitem{struwe.barrett}
M. Struwe, \emph{Wave Maps, in Nonlinear Partial Differential Equations
in Geometry and Physics}, Prog. in Nonlin. Diff. Eq. and their Applic., 
\textbf{29}, (1997), Birkh\"auser, 113--150.

\bibitem{tao:ill}
T. Tao, \emph{Ill-posedness for one-dimensional wave maps at the critical regularity}, Amer. J. Math. 122 (2000), 451--463.

\bibitem{tao:wavemap2}
T. Tao, \emph{Global regularity of wave maps II.  Small energy in two dimensions}, submitted, Comm. Math. Phys.

\bibitem{tataru:wave1}
D. Tataru, \emph{Local and global results for wave maps I}, Comm. PDE \textbf{23} (1998), 1781--1793.

\bibitem{tataru:wave2}
D. Tataru, \emph{On global existence and scattering for the wave maps equation},
Preprint, 1999.

\bibitem{tat:5+1}
D. Tataru, \emph{On $\Box u = |\nabla u|^2$ in $5+1$ dimensions},
Math. Res. Letters \textbf{6} (1999), 469-485.

\end{thebibliography}
\end{document}